\documentclass[12pt,reqno]{article}
\def\hybrid{\topmargin 0pt      \oddsidemargin 0pt
        \headheight 0pt \headsep 0pt
        \textwidth 160true mm       
        \textheight 231true mm         
        \marginparwidth 0.0in
        \parskip 0pt plus 1pt   \jot = 1.5ex}
\usepackage{amssymb}
\usepackage{amsmath}
\usepackage{amsthm}

\hybrid
\usepackage{amssymb}
\usepackage{amsmath}
\usepackage{amsthm}

\newcommand{\be}[1]{\begin{eqnarray#1}}
\newcommand{\ee}[1]{\end{eqnarray#1}}

\newtheorem{thm}{Theorem}[section]

\newtheorem{definition}[thm]{Definition}

\newtheorem{propn}[thm]{Proposition}

\newtheorem{lemma}[thm]{Lemma}

\newtheorem{corollary}[thm]{Corollary}

\renewcommand{\frak}{\mathfrak}

\newcommand{\bd}{\Bar{d}}
\newcommand{\cf}{\hat{f}}
\newcommand{\ot}{\otimes}

\renewcommand{\AA}{\mathcal{A}}
\newcommand{\C}{\mathbb{C}}
\newcommand{\R}{\mathbb{R}}
\newcommand{\Ct}{\mathbb{C}[[t]]}
\newcommand{\Cf}{C_M^\infty}
\newcommand{\CM}{C_M^\infty[[t]]}
\newcommand{\TT}{\mathcal{T}}

\newcommand{\La}{\Lambda}

\newcommand{\ff}{\varphi}
\newcommand{\WW}{\mathcal{W}}
\newcommand{\Wb}{\mathbf{W}}
\newcommand{\EE}{\mathcal{E}}
\newcommand{\PP}{\mathcal{P}}
\newcommand{\QQ}{\mathcal{Q}}

\renewcommand{\Sb}{\mathbf{S}}
\newcommand{\wi}{\mathbf{w}}
\newcommand{\ad}{\mathrm{ad}}
\newcommand{\ft}{\frac{1}{t}}
\newcommand{\OO}{\mathcal{O}}

\newcommand{\g}{\mathfrak{g}}
\newcommand{\rank}{{\rm rank}}
\newcommand{\cl}{{\rm cl}}
\newcommand{\obs}{{{\rm obs}_{DL}}}
\newcommand{\cc}{{\mathbf c}}
\newcommand{\shHom}{\underline{\operatorname{Hom}}}
\renewcommand{\[}{{[\![}}
\renewcommand{\]}{{]\!]}}
\newcommand{\dt}{{\frac{d}{dt}}}

\begin{document}

\title{Polarized deformation quantization}
                          
\author
{Paul Bressler\thanks{Partially supported by the NSF}\\
{\normalsize Max-Planck-Institut f\"ur Mathematik}\\
Joseph Donin\thanks{Partially supported
by Israel Academy of Sciences Grant no. 8007/99 }\\
{\normalsize Dept. of Math.  Bar-Ilan University}\\
{\normalsize Max-Planck-Institut f\"ur Mathematik}}

\date{}
\maketitle
\begin{abstract}
Let $\AA$ be a star product on a symplectic manifold $(M,\omega_0)$,
$\frac{1}{t}[\omega]$ its Fedosov class, where $\omega$ is a deformation of
$\omega_0$. We prove that for a complex polarization of $\omega$
there exists a commutative subalgebra, $\OO$, in $\AA$ that is isomorphic
to the algebra of functions constant along the polarization.

Let $F(\AA)$ consists of elements of $\AA$ whose commutator with $\OO$
belongs to $\OO$. Then, $F(\AA)$ is a Lie algebra which is an $\OO$-extension
of the Lie algebra of derivations of $\OO$. We prove a formula which
relates the class of this extension, the Fedosov class, and the Chern class of $P$.   
\end{abstract}

\section{Introduction}
Let $(M,\omega_0)$ be a symplectic manifold, $\Cf$ the sheaf
of complex valued functions on $M$. A deformation quantization
(or a star-product) on $M$ is a structure of associative algebra on the sheaf of
formal power series $\CM$ with the multiplication of the form
\be{}
f\ast g=\sum_{i\geq 0} t^i\mu_i(f,g), \qquad f,g\in\Cf,
\ee{}
where all $\mu_i$ are bidifferential operators, $\mu_0(f,g)=fg$
and $\mu_1(f,g)-\mu_1(g,f)=\{f,g\}$, the Poisson bracket inverse to $\omega_0$.  

Two deformation quantizations $\AA_1$ and $\AA_2$ on $(M,\omega_0)$
are equivalent if there exists a power series 
$B=1+tB_1+\cdots$, where
$B_i$ are differential operators, such that $B:\AA_1\to\AA_2$ is an
isomorphism of algebras.

It is known that all equivalence classes of star-products on
$M$ with the Poisson bracket $\ff_0=\omega_0^{-1}$
can be obtained by the Fedosov method. According to this method,
one constructs a flat connection, $D$, (called the Fedosov connection) 
on the bundle of Weyl algebras on $M$. 
The quantized algebra, $\AA$, is realized as the subalgebra of
flat sections of the Weyl algebra. The Weyl curvature of $D$, 
being a closed scalar two-form of the view
$\omega=\omega_0+t\omega_1+\cdots$
defines the Fedosov class 
\be{}
\theta(\AA)=\ft[\omega]\in\ft[\omega_0]+H^2(M,\C[[t]]).
\ee{}
It is also known that the correspondence $\AA\mapsto\theta(\AA)$
is a bijection between the set of isomorphism classes of star-products
on $(M,\omega_0)$ and the set $\ft[\omega_0]+H^2(M,\C[[t]])$ modulo
the group of formal symplectomorphisms of $M$, \cite{Fe}, \cite{NT}, \cite{Xu}.

Let $P$ be a polarization of $(M,\omega_0)$. That is,
$P$ is an integrable Lagrangean subbundle of the complexified
tangent bundle $T^\C_M$. Let $\omega=\omega_0+t\omega_1+\cdots$
be a closed form which is a deformation of $\omega_0$ and
$\PP$ a polarization of $\omega$, which is a deformation of polarization $P$.
We say that the triple $(M,\omega,\PP)$ is a deformation of
the triple $(M,\omega_0,P)$.

\begin{definition} 
A {\em polarized star-product}(or quantization) on $(M,\omega_0,P)$ is a pair $(\AA,\OO)$
of $\Ct$ algebras satisfying the conditions:

1) $\AA$ is a star-product on $(M,\omega_0)$, 

2) there exists a deformation $(M,\omega,\PP)$ of $(M,\omega_0,P)$
such that $\OO$ is the sheaf of functions constant along $\PP$,

3) the star-product in $\AA$ being restricted to $\OO$ coincides
with the original commutative multiplication.

Two polarized star-products $(\AA_1,\OO_1)$ and $(\AA_2,\OO_2)$ are
equivalent if there exists an isomorphism $B:\AA_1\to\AA_2$ of star-products
such that $B(\OO_1)=\OO_2$.
\end{definition}

The first result of the paper is the following theorem (cf. Theorem \ref{thmstpr}): 
\begin{thm}\label{th1.2}
For any triple $(M,\omega_0,P)$,
where $P$ is {\em good} (see Definition 2.4) and any its deformation $(M,\omega,\PP)$
there exists a polarized star-product, $(\AA,\OO)$, such that
$\OO$ is the sheaf of functions constant along $\PP$
and $\theta(\AA)=\ft[\omega]$. 
\end{thm}
For example, any K\" ahler or real polarization is good.
This result generalizes a result of N. Reshetikhin and M. Yakimov, \cite{RY},
who consider the case of a real polarization $P$
defined by a Lagrangean fiber bundle $M\to B$. 

Our proof of this result uses the Fedosov method adapted for the case with
polarization. The analogous method
was used by M. Bordemann and S. Waldmann, \cite{BW}, for a construction of
quantization with separation of variables on a K\" ahler manifold.

In fact, we prove a slightly stronger statement. Namely, the 
polarized quantization $(\AA,\OO)$ we constructed satisfies the property
$f\ast g=fg$ for $f\in\OO$ and {\em any} $g\in\AA$.
Such a type of star-products is in the spirit of A. Karabegov
(see \cite{Ka1}, where he constructs star-products with separation of
variables on K\" ahler manifolds).

Let $(\AA,\OO)$ be a polarized quantization on $(M,\omega_0,P)$.
Let $F(\AA)=\{f\in\AA; [f,\OO]\in\OO\}$. The sheaf of $\OO$ modules 
$F(\AA)$ has the
natural structure of a sheaf of Lie algebras with $\OO$ as a commutative
ideal. $T_\OO=F(\AA)/\OO$ can
be considered as a sheaf of Lie algebras consisting of derivations of $\OO$, and 
we may consider the exact sequence of sheaves of $\OO$ modules and Lie algebras
\be{}\label{pls}
0\longrightarrow \OO\longrightarrow F(\AA)\longrightarrow 
T_\OO\longrightarrow 0.    
\ee{}
$F(\AA)$ is an example of an $\OO$-{\em extension of }$T_\OO$, 
(see Definition 4.1)..
We prove that if $P$ is a {em strong} polarization
(see Definition 3.1), this sequence locally splits.
In this case the class $\cl(\AA,\OO)\in H^2(M,\ft\Ct)$ of this extension
is defined.

The second result of the paper is a formula that relates three
classes in $H^2(M,\ft\Ct)$ associated with a polarized deformation
quantization $(\AA,\OO)$: the Fedosov class $\theta(\AA)$,
the Chern class of the polarization $P$,
and $\cl(\AA,\OO)$.
\begin{thm}\label{th1.3}
Let $(\AA,\OO)$ be a polarized quantization of $(M,\omega_0,P)$. Then,
\be{}\label{mF}
\theta(\AA)=\ft\cl(\AA,\OO)-\frac{1}{2}c_1(P).
\ee{}
\end{thm}
One can show that in the case of quantization on a K\" ahler 
manifold with separation of
variables the class $\cl(\AA,\OO)$ coincides with the class defined by Karabegov
in \cite{Ka1}.
The formula analogous to (\ref{mF}) relating the Karabegov and Fedosov
classes was obtained in \cite{Ka2}.

The paper is organized as follows.
In Section 2 we develop a variant of Fedosov's method
adapted to the polarized setting and
prove Theorem \ref{th1.2}.

In Section 3 we define {\em strong} polarizations
and study their local properties.
We show that if $P$ is a strong polarization, then any deformed
polarization, $\PP$, $\PP_0=P$, is also strong, and any strong polarization is
good. We prove that any polarized quantization
of $(M,\omega_0,P)$ with $P$ a strong polarization is locally
isomorphic to a standard quantization, the polarized Moyal star-product. 

In Section 4 we prove Theorem \ref{th1.3}.
To this end we consider, for a given polarized quantization $(\AA,\OO)$, 
the extension (\ref{pls}) and study its behavior under
passing to the opposite quantization, $(\AA^{op},\OO)$,
and to the quantization $(\AA^\sigma,\OO)$ obtained 
by the change of variable $t\mapsto -t$. 
Furthermore, we compare classes $\cl(\AA,\OO)$ with the classes
introduced by Deligne in \cite{De}.

{\bf Acknowledgments.} The authors are grateful to A. Karabegov for
useful discussions and to Max-Planck-Institut f\"ur Mathematik for
hospitality and very stimulating working atmosphere.

\section{Polarized quantization}

\subsection{Setting and notation} 
Let $M$ be a smooth manifold, $C^\infty_M$ the sheaf
of smooth complex valued functions on it. 
Let  $\Ct$ be the algebra of
formal power series over $\C$.
For a sheaf $E$ of $C^\infty_M$ modules, $E[[t]]$ denotes the sheaf
$E\ot_\C \Ct$ completed in the $t$-adic topology.
$E[[t]]$ is a sheaf of $\CM$ modules.

Let $\EE$ be a locally free sheaf of $\CM$ modules of finite rank.
Denote by $T^k(\EE)$ the $k$-th tensor power of $\EE$ over $\CM$
and by $T(\EE)$ the corresponding tensor algebra completed in the
$\{\EE, t\}$-adic topology. 
Similarly we define the completed symmetric algebra $S(\EE)$.
For  a subsheaf $\PP$ of $\EE$ we denote by
$sym_\PP: S(\PP)\to T(\EE)$ the natural map of $\CM$ modules
defined by symmetrization.

Let $\La(\EE)$ be the exterior algebra of $\EE$ over $\CM$.
We will consider $T\ot\La=T(\EE)\ot_{\CM}\La(\EE)$ as a graded super-algebra
regarding a section $x\in T(\EE)\ot\La^k(\EE)$ of degree $k$ even (odd)
if $k$ is even (odd).

Denote by $\delta_T=\delta_{T(\EE)}$ the continuous $\CM$ linear derivation of $T\ot\La$
generated by the map
$T^1(\EE)\ot 1\to 1\ot\La^1(\EE)$, $v\ot 1\mapsto 1\ot v$, $v\in \EE$ is
a section. It is clear that $\delta_T$ is a derivation of degree $1$ and $\delta_T^2=0$.
In the same way define a derivation $\delta_\EE=\delta_{S(\EE)}$ in the (super-)algebra 
$S\ot\La=S(\EE)\ot\La(\EE)$. 
It is easy to see that for a subsheaf $\PP\subset \EE$, 
$\delta_T$ on $T(\EE)\ot\La(\EE)$ can be restricted to 
$S(\PP)\ot\La(\PP)$ via the embedding $sym_\PP\ot id_\PP$,
and coincides with the $\delta_{S(\PP)}$
defined on the algebra $S(\PP)\ot\La(\PP)$.

On the algebra $S\ot\La$ there is defined another
derivation, $\delta^*$, of degree $-1$ generated by the map $1\ot\La^1\to S^1\ot 1$,
$1\ot v\to v\ot 1$. It is easy to check that
$(\delta^*)^2=0$ and $[\delta,\delta^*]=\delta\delta^*+\delta^*\delta=deg$, where
$deg$ is the derivation assigning to an element $x\in S^p\ot\La^q$ the
element $(p+q)x$. 

Let $\EE$ be presented as a direct sum of $\CM$ submodules,
$\EE=\PP\oplus\QQ$. We will identify the tensor product of algebras,
$(S(\PP)\ot\La(\PP))\ot (S(\QQ)\ot\La(\QQ))$, with
$(S(\PP)\ot S(\QQ))\ot (\La(\PP)\ot \La(\QQ))=(S(\PP)\ot S(\QQ))\ot \La(\EE)$,
where the last identification is obtained via skew-symmetrization of the last factor.
Denote $S(\PP,\QQ)=S(\PP)\ot S(\QQ)$, $\Sb=\Sb(\PP,\QQ)=(S(\PP)\ot S(\QQ))\ot \La(\EE)$, 
$\Sb(\PP)=S(\PP)\ot 1\ot \La(\EE)$,
and $\Sb(\QQ)=1\ot S(\QQ)\ot\La(\EE)$.

One has the embedding 
\be{}\label{rels}
sym_\PP\ot sym_\QQ\ot id: (S(\PP)\ot S(\QQ))\ot \La(\EE)\to T(\EE)\ot\La(\EE).
\ee{}

It is obvious that the derivations $\delta_\PP$,
$\delta^*_\PP$, $\delta_\QQ$, $\delta^*_\QQ$ induce
the corresponding derivations on the algebra $\Sb$,
and $\delta=\delta_\PP+\delta_\QQ$ coincides with
the restriction to $\Sb$ of $\delta_T$ on $T(\EE)\ot\La(\EE)$ via the embedding (\ref{rels}). 

Let us define the operator $\delta^{-1}=\delta^{-1}_{\PP,\QQ}$ on $\Sb$.
There is the continuous $\CM$ linear map that is equal to zero for $x\in \CM$
and $\delta^{-1}(x)=(1/(p+r+q)\delta^*(x)$ 
for $x\in  (S^p(\PP)\ot S^r(\QQ))\ot \La^q(\EE)$, $p+r+q>0$.

There is the obvious relation
\be{}\label{relde}
\delta\,\delta^{-1}+ \delta^{-1}\delta= \mbox{ projection on $\Sb^+$ along $\CM$},
\ee{}
where $\Sb^+$ is the closure of $\oplus_{p+r+q>0}(S^p(\PP)\ot S^r(\QQ))\ot \La^q(\EE)$.

\subsection{Fedosov algebra}
Let $\ff:\EE\ot\EE\to\CM$ be a $\CM$ linear skew-symmetric form
and $I$ the closed ideal in $T(\EE)$ generated by relations
\be{}\label{relI}
x\ot y-y\ot x=t\ff(x,y).
\ee{}
We call $\WW(\EE)=T(\EE)/I$ the Weyl algebra and
$\Wb=\Wb(\EE)=\WW\ot\La(\EE)$ the Fedosov algebra over $\EE$.

Let $\EE=\PP\oplus\QQ$ be a decomposition into $\CM$ modules.

The derivation $\delta$ of $T(\EE)\ot\La(\EE)$ induces
a derivation of $\Wb$. Indeed, $\delta$ applied to
the both sides of (\ref{relI}) gives zero.

Define the Wick map, $\wi=\wi_{\PP,\QQ}$, as the composition 
$\Sb(\PP,\QQ)\to T(\EE)\ot\La(\EE)\to \Wb$, where the first map
is (\ref{rels}) and the 
second is the projection. By the PBW theorem 
$\wi$ is an isomorphism of $\CM$ modules.

Due to the isomorphism $\wi$, all the operators
$\delta_\PP$, $\delta_\QQ$, $\delta^{-1}_{\PP,\QQ}$
carry over from $\Sb$ to $\Wb$.
We retain for them the same notation. 
Note that while $\delta_\PP+\delta_\QQ$
does not depend on
the decomposition and coincides with the $\delta$ induced from
$T(\EE)\ot\La(\EE)$, $\delta^{-1}_{\PP,\QQ}$ is not a derivation
and does depend on the decomposition. 
In particular, one can suppose that the decomposition is trivial,
$\EE=\EE\oplus 0$. In this case we denote
$\delta^{-1}_\EE=\delta^{-1}_{\EE,0}$.

Note that $\delta_\PP$, $\delta_\QQ$, and $\delta$ are derivations on $\Wb$.

\begin{propn}\label{prop1.1} One has
$$H(\Wb,\delta)=\CM.$$
Moreover, if $x\in \WW(\EE)\ot\La^{k>0}(\EE)$
then $y=\delta^{-1}_{\PP,\QQ}x$ is such that
$\delta y=x$ for any decomposition $\EE=\PP\oplus\QQ$.
\end{propn}
\begin{proof} Follows from (\ref{relde}).
\end{proof}

\subsection{Lie subalgebras in $\WW$}
Let $\EE=\PP\oplus\QQ$ be a decomposition.
We say that $x\in \Wb$ has $\wi_{\PP,\QQ}$-degree $(p,q)$ if
$\wi_{\PP,\QQ}^{-1}(x)\in (S^p(\PP\ot S^q(\QQ))\ot\La(\EE)$.
We say that $x\in \Wb$ has $\wi_{\PP,\QQ}$-degree $k$ if
$\wi_{\PP,\QQ}^{-1}(x)\in \oplus_{p+q=k}(S^p(\PP)\ot S^q(\QQ))\ot\La(\EE)$.
The $\wi_\EE$-degree is $\wi_{\EE,0}$-degree for the trivial
decomposition $\EE=\EE\oplus 0$.

Let $\ff$ be nondegenerate.
Let $\g$ be a sheaf of Lie algebras acting on $\EE$.
We call a $\CM$ linear map $\lambda:\g\to\WW$ a realization of
$\g$ if it is a Lie algebra morphism ($\WW$ is considered as a Lie algebra with respect
to commutator $\ft [\cdot,\cdot]$) and for any $x\in\g$ and $e\in\EE$ 
one has $x(e)=\ft[\lambda(x),e]$. Any two realizations differ
by a Lie algebra morphism of $\g$ to the center of $\WW$, so
if $\g$ is a sheaf of semisimple Lie algebras 
there is not more than one realization of $\g$.

Denote by ${\frak{sp}}(\EE)$ the sheaf of symplectic Lie algebras
with respect to $\ff$. Since ${\frak{sp}}(\EE)$ is semisimple, there is a unique
realization  $\rho_\EE:{\frak{sp}}(\EE) \to \WW$. The image of this realization 
consists of elements having $\wi_\EE$-degree two.
 
Let $\EE=\PP\oplus\QQ$ be a decomposition into Lagrangean subsheaves.
Denote by ${\frak{sp}}(\PP,\EE)$ the subsheaf of ${\frak{sp}}(\EE)$
preserving $\PP$. 
It is easy to check that ${\frak{sp}}(\PP,\EE)$ can be realized as the subset
of elements of $\WW$ having $\wi_{\PP,\QQ}$-degree $(1,1)$ and $(2,0)$.
Denote this realization by $\rho_{\PP,\EE}:{\mathfrak{sp}}(\PP,\EE)\to \WW$.
On the other hand, ${\mathfrak{sp}}(\PP,\EE)$ can be realized in $\WW$ by $\rho_\EE$.

Let us define a $\C[[t]]$ linear map $tr: {\mathfrak{sp}}(\PP,\EE)\to \C[[t]]$
in the following way.
Let $a\in {\frak{sp}}(\PP,\EE)$ 
and $a'$ its restriction to $\PP$. We put $tr(a)=trace(a')$.
It is easy to see that $\rho_{\PP,\EE}(a')$  (here $a'$ is trivially extended 
to $\EE$) is the $(1,1)$ component of $\rho_{\PP,\EE}(a)$. 

\begin{lemma}\label{lem1.1}
Let $\EE=\PP\oplus\QQ$ be a decomposition into Lagrangean subsheaves.
Let $a\in {\frak{sp}}(\PP,\EE)$.
Then 
$$\rho_{\PP,\EE}(a)-t\frac{1}{2}tr(a)=\rho_\EE(a).$$

\end{lemma}

\begin{proof} Straightforward.
\end{proof}

\subsection{Filtrations on $\WW$}
We define two decreasing filtrations on $\WW$
numbered by nonnegative integers.

The $T$-filtration $F^T_\bullet\WW$ is defined as follows.
We ascribe to the elements of $\EE$ degree 1 and to $t$ degree 2.
Then $F^T_n\WW$  consists of elements of $\WW$ 
having the leading term of total degree $\geq n$.

The $\PP$-filtration, $F^\PP_\bullet\WW$, is firstly defined on $S(\PP)\ot S(\QQ)$
by the subsets
$F^\PP_n=S^n(\PP)S(\PP)\ot S(\QQ)$, $n=0,1,...$, and carried over to $\WW$ via the Wick
isomorphism.  

We extend those filtrations to $\Wb$ in the natural way standing,
for example, $F^T_n\Wb=F^T_n\WW\ot\La(\EE)$.
We will use the following mnemonic notation. To point out, for example,
that a section $x\in \Wb$ belongs to $F^T_n\Wb$ we write
$F^T(x)\ge n$.

 We call a subsheaf $\PP\subset\EE$ a $\ff$-null subsheaf
if $\ff$ when restricted to $\PP$ equals zero.

\begin{propn}\label{prop1.2} 
Let $\EE=\PP\oplus \QQ$ be a decomposition of $\EE$ with $\PP$ a
$\ff$-null subsheaf. Then

a) The Wick map $\wi:\Sb\to\Wb$ has the following property:
for $a\in \Sb(\PP)$
and arbitrary $c\in\Sb$ one has
$ac=\wi(a)\wi(c)$.

 The filtrations on $\Wb$ have the properties:

b) for $x,y\in \Wb$, if $F^\PP(x)\ge k$ then $F^\PP(xy)\ge k$;

c) $F^\PP(\delta^{-1}_{\PP,\QQ}x)\geq F^\PP(x)$;

d) $F^T(\delta^{-1}_{\PP,\QQ}x)>F^T(x)$.
\end{propn}

\begin{proof} a) follows from the fact that $\PP$ is a Lagrangean
subsheaf of $\EE$ and the definition of the Wick map.
b) follows from a), c) and d) are obvious.
\end{proof}

\subsection{Good polarizations and connections}

Let $T_M$ be the complexified tangent bundle over $M$, $\TT=T_M[[t]]$, 
and $\omega:\TT^{\ot 2}\to\CM$ a nondegenerate symplectic form.
It means that $\omega=\omega_0+t\omega_1+\cdots$
where $\omega_0$ is nondegenerate and $d\omega=0$. 
Let $\ff=\ff_0+t\ff_1+\cdots$ be 
the nondegenerate Poisson bracket
inverse to $\omega$. It may be considered as the map
$\ff:(\TT^*)^{\ot 2}\to\CM$, where $\TT^*=T_M^*[[t]]$. 

For any function $f\in\CM$ denote by $X_f$ the corresponding
Hamiltonian vector  field on $M$, $X_f=\ff(df,\cdot)$.
The map $\phi:\TT^*\to\TT$ defined by $df\mapsto X_f$ is
an isomorphism of sheaves. Its inverse map is given by
$X_f\mapsto\omega(\cdot,X_f)=df$.

\begin{definition}\label{def1.1}
a) A direct subsheaf $\PP\subset\TT$ is called Lagrangean
if $\omega(\PP,\PP)=0$ and $rank(\PP)=\frac{1}{2}\dim M=n$.

b) $\PP$ is called integrable if locally there exist functions
$a_1,...,a_n$ such that
the Hamiltonian vector fields $X_{a_1},...,X_{a_k}$
form a local basis in $\PP$.  

c) A Lagrangean integrable $\PP$ is called a polarization of $\omega$.
In this case all the $X_{a_i}$, $i=1,...,n$, pairwise commute,
and $da_i$ form a local basis in $\PP^\perp$.

d) A polarization $\PP$ is called good if
there exist additional functions $f_1,...,f_n$ such that
the Hamiltonian vector fields
$X_{a_i}$, $X_{f_i}$, $i=1,...,n$, pairwise commute
and form a local basis in $\TT$.
\end{definition} 

\begin{definition}
Let $\PP\subset\TT$ be a Lagrangean subsheaf.
We call a connection, $\nabla$, on $M$ 
a $\PP$-{\em connection} if

a) it preserves $\omega$ and is torsion free, i.e.
is a symplectic connection, and

b) it preserves $\PP$, i.e. $\nabla(\PP)\subset\PP\ot\La^2(\TT^*)$,
and is flat on $\PP$ along $\PP$, i.e.
for any $X,Y\in\PP$ one has
$(\nabla_X\nabla_Y-\nabla_Y\nabla_X-\nabla_{[X,Y]})(\PP)=0$.
\end{definition}

\begin{propn}\label{prop1.6}
Let $\PP\subset\TT$ be a good polarization of $\omega$.
Then, there exists a $\PP$-connection on $M$. 
\end{propn}

\begin{proof} 
Let $a_1,...,a_{2n}$ be functions such that $X_i=X_{a_i}$, $i=1,...,n$, form
a local basis in $\PP$ and all the $X_i$ commute and form a local basis in $\TT$.
Let us set $\nabla_{X_i}X_j=[X_i,X_j]=0$. This rule defines a local connection.

Let us prove the invariance of $\ff$, i.e.
$$\nabla_z(\omega(x,y))=\omega(\nabla_z x,y)+\omega(x,\nabla_z y) \mbox{  for any }x,y,z\in \TT.$$
But for $x=X_i$, $y=X_j$, $z=X_k$ this equation
is equivalent to the relation 
$$\ff(a_k,\ff(a_i,a_j))=\ff(\ff(a_k,a_i),a_j)+\ff(a_i,\ff(a_k,a_j))=0,$$
which holds by the Jacobi identity for the Poisson bracket $\ff$.

Let us proof that $\nabla$ is torsion free.
Since the torsion is $\CM$ linear, it is enough to prove that
$$\nabla_{X_i}(Y_j)-\nabla_{X_j}(X_i)-[X_i,X_j]=0$$
for all pairs $X_i,X_j$ of our basis.
But this follows from the definition of $\nabla$ and from the fact
that all $X_i$ pairwise commute.
Similarly, from the fact
that all $X_i$ pairwise commute follows that $\nabla$ is flat.

That $\nabla$ preserves $\PP$ is obvious.

Since $X_f\in\PP$ is equivalent to
$df\in\PP^\perp$, it is obvious that
 $\nabla$ has the property: for any $f,g\in\CM$ such that
$X_f,X_g\in \PP$ one has
\be{}\label{hamilt}
\nabla_{X_f} X_g=0.
\ee{}

Now, let us prove the existence of a global connection.

Let $\{U_\alpha\}$ is an open covering of $M$ such that
on each $U_\alpha$ there is a $\PP$-connection $\nabla_\alpha$.
Then the differences $\nabla_\alpha-\nabla_\beta$ form on
$U_\alpha\cap U_\beta$ a \^ Cech cocycle 
$\psi_{\alpha,\beta}\in Hom(\TT\ot\TT, \TT)$,
$\psi_{\alpha,\beta}(X,Y)=\nabla_{\alpha X} Y-\nabla_{\beta X} Y$.

$\psi_{\alpha,\beta}$ satisfy the following properties.

It follows from (\ref{hamilt}) that
\be{}\label{hamilt1}
\psi_{\alpha,\beta}(X,Y)=0 \quad\quad \mbox{for\ } X,Y\in\PP.
\ee{}
Since all $\nabla_\alpha$ are torsion free, $\psi_{\alpha,\beta}$
are symmetric. Since all $\nabla_\alpha$ preserve $\PP$,
one has $\psi_{\alpha,\beta}(X,Y)\in \PP$ for $Y\in\PP$.

In addition, $\psi_{\alpha,\beta}$ considered as elements from
$Hom(\TT, Hom(\TT,\TT)$, $X\mapsto\psi_{\alpha,\beta}(X,\cdot)$,
belong to $Hom(\TT,sp(\TT))$ where $sp(\TT)$ consists
of endomorphisms of $\TT$ preserving $\omega$ invariant.

Since all the properties above are $\CM$ linear, it is possible to find
tensors $\psi_\alpha\in Hom(\TT\ot\TT, \TT)$ satisfying the
same properties and such that $\psi_\alpha-\psi_\beta=\psi_{\alpha,\beta}$.
Then $\nabla=\nabla_\alpha-\psi_\alpha=\nabla_\beta-\psi_\beta$
is a globally defined connection.
Flatness of $\nabla$ on $\PP$ along $\PP$ follows from
the fact that for all $\alpha$ $\psi_\alpha(X,Y)=0$ for $X,Y\in\PP$, which follows from 
(\ref{hamilt1}). $\nabla$ is torsion free because all $\psi_\alpha$ are symmetric.
So, $\nabla$ satisfies the proposition.
\end{proof}

\subsection{Fedosov's construction}
\label{FC}
Let $(M,\omega_0)$ be a symplectic manifold. 
It is known that all equivalent classes of star-products on
$M$ with the Poisson bracket $\ff_0=\omega_0^{-1}$
can be obtained by the Fedosov method. According to this method,
one constructs a flat connection, $D$, (called the Fedosov connection) 
in the Weyl algebra defined on 
the cotangent bundle with the help of relations (\ref{relI}),
where $\ff=\ff_0$. The quantized algebra, $\AA$, is realized as the subalgebra of
flat sections of the Weyl algebra. The Weyl curvature of $D$, 
being a closed scalar two-form of the view
$\tilde\theta=\omega_0+t\omega_1+\cdots$
defines the Fedosov class 
\be{}\label{formtheta}
\theta(\AA)=\ft[\tilde\theta]\in\ft[\omega_0]+H^2(M,\C[[t]]).
\ee{}
It is also known that the correspondence $\AA\mapsto[\theta(\AA)]$
is a bijection between the set of isomorphism classes of star-products
on $(M,\omega_0)$ and the set $\ft[\omega_0]+H^2(M,\C[[t]])$ modulo
the group of formal symplectomorphisms of $M$, \cite{Fe}, \cite{NT}, \cite{Xu}.

We will adapt the Fedosov method to construct a polarized star-product
corresponding to a class $\ft[\tilde\theta]\in\ft[\omega_0]+H^2(M,\C[[t]])$
with a polarization of $\tilde\theta$.
We will see that in the presence of polarization by realizing the Fedosov scheme 
there
appears also a Wick curvature which differs from the Weyl
curvature by $t$ times a half of the first Chern class of the polarization. 

Let $\PP'\subset\TT$ be a good polarization of $\omega$.
It is easy to prove that there exists a complement Lagrangean subsheaf 
$\QQ'\subset\TT$ such that $ \TT=\PP'\oplus\QQ'$ on $M$.

In the following we set $\EE=\TT^*$ and consider
the decomposition $\EE=\PP\oplus\QQ$, there $\PP$ and $\QQ$
correspond to $\PP'$ and $\QQ'$ by the isomorphism 
$\phi$ between $\TT$ and $\EE$. Note that
$\PP$ and $\PP'$ are orthogonal to each other, so that we
may write $\PP'=\PP^\perp$.

Let $\nabla$ be a $\PP^\perp$-connection on $M$.
Then the induced connection $\nabla:\EE\to\EE\ot\La^1\EE$
on $\EE$ preserves $\PP$, i.e.
 $\nabla(\PP)\subset\PP\ot\La^2(\EE)$,
and is flat on $\PP$ along $\PP'$, i.e.
for any $X,Y\in\PP'$ one has
$(\nabla_X\nabla_Y-\nabla_Y\nabla_X-\nabla_{[X,Y]})(\PP)=0$.

$\nabla$ gives a derivation of the Fedosov algebra $\Wb=\Wb(\EE)$, 
which is an extension
of the de Rham differential on functions. 
Analogously, $\nabla$ gives such derivations of the algebras $T(\EE)\ot\La(\EE)$,
$S(\EE)\ot\La(\EE)$, and $(S(\PP)\ot S(\QQ))\ot\La(\EE)$.
These derivations commutes with the maps
(\ref{rels}) and $\wi$.

For convenience, we will mark the elements of the Fedosov algebra
lying in $\EE\ot 1$ by letters with hat over them ($\hat{x}$), and by $dx$ we will
denote the copy of $\hat{x}$ lying in $1\ot\La^1 \EE$.

It is easy to check that for 
$\tilde{\delta}=\omega_{ij}\Hat{x^i}\ot dx^j$ one has
\be{}\label{reldelta}
\delta=\frac{1}{t}\ad(\tilde{\delta}),  \notag \\   
\tilde{\delta}^2=t\omega.
\ee{}
From the fact that the torsion of $\nabla$ is equal to zero follows
\be{}\label{torsion}
\nabla(\tilde\delta)=0.
\ee{}

Since $\nabla^2$ is a $\CM$ linear derivation of degree $0$ preserving $\PP$, 
there is an element $R\in \rho_{\PP,\EE}(\frak{sp})\ot\La^2(\EE)$ 
 such that 
$\nabla^2=\frac{1}{t}\ad(R)$. In particular,
one has to be 
\be{}\label{FPR}
F^\PP(R)\geq 1.
\ee{}

According to Fedosov, \cite{Fe}, we also define $R^F\in\rho_\EE({\frak{sp}})\ot\La^2(\EE)$
satisfying $\nabla^2=\ft\ad(R^F)$.

From (\ref{torsion}) follows
\be{}
\delta(R)=\delta(R^F)=0.
\ee{}

Following to Fedosov, we will consider connections on $\Wb$ of the form
\be{}
D=\nabla+\ft\ad(\gamma), \quad\quad \gamma\in \WW\ot\La^1(\EE).
\ee{}
We define the Wick curvature of $D$ as
\be{}\label{wickcur}
\Omega_D=R+\nabla(\gamma)+\ft\gamma^2. 
\ee{}

According to Fedosov, we also define the Weyl (or Fedosov) curvature of $D$ as
\be{}\label{weylcur}
\Omega_D^F=R^F+\nabla(\gamma)+\ft\gamma^2=
R-t\frac{1}{2}tr(R)+\nabla(\gamma)+\ft\gamma^2=\Omega_D-t\frac{1}{2}tr(R), 
\ee{}
where the second equality is due to Lemma \ref{lem1.1}.

\begin{propn}\label{propn1.2}
tr(R) is a closed form polarized by $\PP$. The class of $tr(R)$ 
in $H^2(M,\C)$ coincides with
the first Chern class of $\PP$, i.e. one has
\be{}\label{chcl} 
[tr(R)]=c_1(\PP)=c_1(\PP_{t=0}).
\ee{}
\end{propn}

\begin{proof}
Since $tr(R)$ is a $\CM$ valued 2-form, it follows
from the Bianchi identity that it is a closed 2-form.
From flatness of $\PP$ along $\PP^\perp$ follows that
the $(1,1)$ component of $R$ does not contain terms of
the view $R_{i,j}\ot dy^i\wedge dy^j$ with $dy^i,dy^j\in \QQ$. 
It follows that $\PP$ is a polarization of $tr(R)$.

It is easy to see that
$\nabla$ being restricted to $\PP$ defines a connection on 
$\PP$, $\nabla_\PP$,
and 
$$\nabla^2_\PP=\ft \ad R_{1,1},$$
where $R_{1,1}$ is the $\wi_{\PP,\QQ}$-degree $(1,1)$ 
component of the $\frak{sp}_{\PP,\QQ}$ valued 
2-form $R$.
Now (\ref{chcl}) follows from the definition of the Chern
classes. Since Chern classes are integer valued elements of
$H^2(M,\C)$, i.e. belong to $H^2(M,\mathbb{Z})$, the
independence of $c_1(\PP)$ on the parameter $t$ is obvious.
\end{proof}
   
One has
\be{}\label{flatness}
D^2=\ft\ad(\Omega_D)=\ft\ad(\Omega_D^F).
\ee{}

Let us take $\gamma$ in the form 
\be{}
\gamma=\tilde\delta+r, \quad\quad r\in \WW\ot\La^1(\EE), \quad F^T(r)\geq 3.
\ee{}
Then the connection $D$ has the form
\be{}\label{Fcon}
\nabla+\delta+\ft\ad(r).
\ee{}
Using (\ref{reldelta}) and (\ref{torsion}), we obtain that its
Wick curvature is
\be{}\label{Fcurv}
\Omega_D=R+\nabla(\tilde\delta+r)+\frac{1}{t}(\tilde\delta+r)^2=
\omega+\delta r+R+\nabla r+\ft r^2.
\ee{}

\begin{propn}\label{prop1.3}
There exists an element $r\in \WW(\EE)\ot\La^1(\EE)$ such that

a) $F^T(r)\geq 3$;

b) $F^\PP(r)\geq 1$;

c) the connection $D=\nabla+\delta+\ft\ad(r)$ is flat, i.e.
$D^2=0$;

d) for its Wick curvature one has 
$$\Omega_D=\omega;$$

e) for the class in $H^2(M,\Ct)$ of its Weyl curvature one has 
\be{}\label{Wc}
[\Omega_D^F]=[\Omega_D]-t\frac{1}{2}c_1(\PP).
\ee{}
\end{propn}

\begin{proof}
First of all, we apply the Fedosov method, (\cite{Fe}, Theorem 5.2.2),
to find $r$ satisfying d).  
According to (\ref{Fcurv}), $r$ must obey  the
equation
\be{}
\delta r=-(R+\nabla r+\ft r^2)
\ee{}
Look for $r$ as the limit of the sequence, $r=\lim r_k$, 
where $r_k\in \WW(\EE)\ot\La^1(\EE)$, $k=3,4,...$,
and $F^T(r_{k}-r_{k-1})\geq k$.
As in Lemma 5.2.3 of \cite{Fe}, 
using Proposition (\ref{prop1.2}) d) and the fact that $F^T(R)\geq 2$,
such $r_k$ can be calculated recursively:
\be{}\label{calr}
r_3=-\delta_{\PP,\QQ}^{-1}(R) \notag\\
r_{k+3}=-(r_3+\delta_{\PP,\QQ}^{-1}(\nabla  r_{k+2}+\frac{1}{t}r_{k+2}^2)).
\ee{}
So, a) and d) are proven. e) follows from (\ref{weylcur}) and Proposition \ref{propn1.2}.

Let us prove that $F^\PP(r_k)\geq 1$ for all $k$.
That $F^\PP(r_3)\geq 1$  follows from the fact
that $F^\PP(R)\geq 1$ and from Proposition (\ref{prop1.2}) c).
Suppose that $F^\PP(r_i)\geq 1$ for $i<k+3$. Then 
$F^\PP(\nabla r_{k+2})\geq 1$ because $\nabla$ preserves $\PP$.
On the other hand, $F^\PP(r_{k+2}^2)\geq 1$
because of Proposition (\ref{prop1.2}) b),
therefore from (\ref{calr}) follows that $F^\PP(r_{k+3})\geq 1$ as well. 
So, we have that $r$ being the limit of the convergent sequence
$r_k$ satisfy the conditions a), b), and d) of
the proposition. c) obviously follows from d) and  (\ref{flatness}).
\end{proof}

Denote by $\WW_D$ the subsheaf of $\WW$ consisting of
flat sections $a$, i.e. such that $Da=0$.
Since $D$ is a derivation of $\Wb$, it is clear that $\WW_D$ is a
sheaf of subalgebras. 
Let $\sigma=id-(\delta\delta_{\PP,\QQ}^{-1}+\delta_{\PP,\QQ}^{-1}\delta)$.
Then, as follows from (\ref{relde}), $\sigma: \WW\to\CM$ is a projection,
where $\CM$ considers as the center of the algebra $\WW$. 

\begin{propn}\label{prop1.4}
a) The map $\sigma:\WW_D\to\CM$ is a bijection.

b) The inverse map $\tau:\CM\to \WW_D$ has the form
$\tau(f)=f+\hat f$, there $F^T(\hat f)>F^T(f)$.

c) If $df\in\PP$ then $F^\PP(\hat f)\geq 1$.

d) If $df\in\PP$ then $\sigma(\tau(f)\tau(g))=fg$
for any $g\in\CM$.
\end{propn}

\begin{proof}
Again, we apply the Fedosov iteration procedure.
According to \cite{Fe}, Theorem 5.2.4, we look for $\tau(f)$
as a limit, $\tau(f)=\lim a_k$, there $a_k\in \WW$ can
be calculated recursively:
\be{}\label{calf}
a_0=f \notag \\ 
a_{k+1}=a_0+\delta_{\PP,\QQ}^{-1}(\nabla a_k+\ft\ad r(a_k)).
\ee{}
Put $\hat f=\tau(f)-a_0$.
As in \cite{Fe}, Theorem 5.2.4, one proves that such $\tau(f)$ 
and $\hat f$ satisfy a) and b).
Now observe that $a_1-a_0=\delta_{\PP,\QQ}^{-1}(1\ot df)$ and
if  $df\in\PP$ then $F^\PP(a_1-a_0)\geq 1$.
By induction we conclude that  $F^\PP(a_k-a_0)\geq 1$ for all $k\geq 1$.
So $F^\PP(a-a_0)\geq 1$ as well, which proves c).

Let us prove d). We have $\tau(f)\tau(g)=f\tau(g)+\hat f\tau(g)$.
Since by c) $F^\PP(\hat f)\geq 1$, $F^\PP(\hat f\tau(g))\geq 1$ as well.
It follows that $\sigma(\hat f\tau(g))=0$ and 
$\sigma(\tau(f)\tau(g))=\sigma(f\tau(g))=fg$,
because $\sigma$ is a $\CM$ linear map and $\sigma(\tau(g))=g$. 
\end{proof}

\subsection{Polarized star-product}
Let $T_M$ be the complexified tangent bundle over $M$, $\TT=T_M[[t]]$, 
and $\omega:\TT^{\ot 2}\to\CM$ a nondegenerate symplectic form,
$\omega=\omega_0+t\omega_1+\cdots$.
Let $\PP\subset\TT$ be a polarization of $\omega$.
We say that $\omega$ is a deformation of the symplectic form $\omega_0$
and $\PP$ is a deformation (or extension) of $\PP_0=\PP/t\PP$.
It is clear that $\PP_0$ is a polarization of $\omega_0$.

Let $\tilde\omega$ be another deformation of $\omega_0$.
We say that $\omega$ and $\tilde\omega$ are equivalent
if they define the same class in $\omega_0+H^2(M,t\C[[t]])$. 

Let $\PP$ be a polarization of $\omega$ which is an extension of $\PP_0$.
It is easy to show that if $\Tilde\omega=\omega+d\lambda$ and
$X_\lambda$ is a vector field corresponding to $\lambda$ then
$\Tilde\PP=\exp(L_{X_\lambda})\PP$ is a polarization of $\Tilde\omega$
extending $\PP_0$, which we call an equivalent polarization.
Thus, it is meaningful to speak about a polarization of a
class of $\omega_0+H^2(M,t\Ct)$. It is the equivalence class
of extensions of $\PP_0$. 

As we have mentioned at the beginning  of Section \ref{FC}, 
the classification theorem for star-products say that
all star-products on $\CM$ with the Poisson bracket inverse to $\omega_0$
are parameterized by elements 
\be{}
\theta=\ft[\tilde\theta]\in \ft[\omega_0]+H^2(M,\Ct)
\ee{}
where $\tilde\theta$ is the Weyl scalar curvatures of the Fedosov connections.

\begin{thm}\label{thmstpr}
Let $(M,\omega_0,)$ be a symplectic manifold.
Let $\omega$ be a deformation of $\omega_0$
and $\PP$ a good polarization of $\omega$.
Let $\OO$ be the sheaf of functions from $\CM$
constant along $\PP$. 

Then, there exists a star-product, $\ast$, on $\CM$
having the Fedosov class $\ft[\omega]$ and satisfying the property:
\be{}\label{psp}
f\ast g=fg \qquad \mbox{for all } f\in\OO.
\ee{}
\end{thm}

\begin{proof}
Let $\alpha$ be a 2-form representing  $\frac{1}{2}c_1(\PP)$ 
and having $\PP$ as its polarization 
(see Proposition (\ref{propn1.2})).
Construct the star-product for 
$\omega'=\omega+t\alpha$ and $\PP$
as in Proposition \ref{prop1.4}.
According to Proposition \ref{prop1.3} d), e), the Wick
curvature corresponding to that star-product is equal to $\omega'$
and its Fedosov class is equal to $\ft[\omega]$. 

Property (\ref{psp}) follows from Proposition \ref{prop1.4} d).
\end{proof}

\section{Strong polarizations}

\subsection{Definition and deformation of a strong polarization}
Let $(M,\omega_0)$ be a symplectic manifold of dimension $2n$, 
$T$ the tangent bundle on $M$ and $T\ot_\R\C=T^\C$
its complexification. 
Let $P\subset T^\C$ be a polarization of $\omega_0$, i.e.
$P$ is a Lagrangean subbundle and locally there exist
functions $a_1,...,a_n$ such that $da_1,...,da_n$
give a local basis for $P^\perp\subset (T^\C)^*$.
Let the triple $(M,\omega,\PP)$ be a deformation of the triple $(M,\omega_0,P)$,
i.e. $\PP$ is a polarization of symplectic form $\omega=\omega_0+t\omega_1+\cdots$, 
$\PP_0=P$.
We are going to give simple conditions on $P$ and $\PP$
which guarantee for them to be good polarizations.

Denote by $\OO=\OO_P\subset\Cf$ the subsheaf of functions constant
along the polarization $P$.
In the obvious way one can form 
an analog of
the Dolbeault complex, which we also call Dolbeault:
$(\bigwedge^\bullet P^*, \Bar d)$, where $\Bar d$ is the differential
along $P$.
By definition, $\OO$ consists of functions $f$ such that $\Bar df=0$.
Note that the Dolbeault complex is meaningful
for any involutive $P$, i.e. when $[P,P]\subset P$, and for
any involutive defirmation $\PP$.

\begin{definition}\label{defn3.1} We call a polarization $P$ 
{\em strong} if for any point $m\in M$ there exists a neighborhood $U$ of $m$
such that
the complex of sections $(\Gamma(U;\bigwedge^\bullet P^*), \Bar d)$
is exact and gives a resolution of $\Gamma(U;\OO_P)$.

We give the same definition  for a deformed polarization $\PP$
of $\omega$.
\end{definition}

\begin{propn}\label{prop3.1}
Let $P$ be a strong polarization of $\omega_0$.
Let $\omega$ be a symplectic
deformation of $\omega_0$
and $\PP\subset\TT=T^\C[[t]]$ a deformation of $P$, $\PP_0=P$,
which is  Lagrangean  with respect to $\omega$
and involutive, i.e. $[\PP,\PP]\subset\PP$.
Then $\PP$ is a strong polarization of $\omega$.
\end{propn}

\begin{proof}
First of all, let us prove that $\PP$ is a polarization, i.e. $\PP$ is integrable.
Since $\PP$ is involutive, the Dolbeault complex $(\bigwedge^\bullet \PP^*, \Bar d_t)$
is meaningful, where $\Bar d_t=\Bar d+t\Bar d'$ is a deformation of
the differential $\Bar d$ from the Dolbeault complex for $P$.
Let $U$ be a neighborhood such that the complex
$(\Gamma_U(\bigwedge^\bullet P^*), \Bar d)$ is exact.
It is enough to prove that if $f$ is a function on $U$ satisfying $\bd f=0$ then
there is an extension $f_t=f+tf'\in C^\infty_U[[t]]$ such that $\bd_tf_t=0$.

Suppose that we have already found a series $f^{(n-1)}=f+tf_1+\cdots+t^{n-1}f_{n-1}$
such that $\bd_tf^{(n-1)}=0$ mod $t^{n}$.
Then $\bd_tf^{(n-1)}=t^n\gamma_n$ mod $t^{n+1}$.
Applying $\bd_t$ to the both sides of this equation and dividing by $t^n$ we obtain
$\bd_t\gamma_n=0$ mod $t$. Hence, $\bd\gamma_n=0$.
So, $\gamma_n$ is a closed 1-form in the Dolbeault complex corresponding to $P$,
and one can find
a function $f_n$ such that $\bd f_n=-\gamma_n$.
It is clear that if we put $f^{(n)}=f^{(n-1)}+t^nf_n$,
we obtain $\bd_tf^{(n)}=0$ mod $t^{n+1}$.
In this way we construct step-by-step the series
$f_t$ such that $\bd_tf_t=0$ and prove
that $\PP$ is a polarization of $\omega$. 

That $\PP$ is strong is obvious because of the upper-semicontinuity
argument.
\end{proof}

Proposition \ref{prop3.1} shows that the deformed polarization $\PP$ is strong
if $P=\PP_0$ is strong.
There is the following proposition provides sufficient
conditions
for a Lagrangean subbundle $P$ to be a strong polarization.

\begin{propn} 
Let $(M,\omega_0)$ be symplectic manifold with a Lagrangean subbundle
$P\subset T_M^\C$. Suppose $P$ satisfies the following conditions:

i) the subsheaf $P\cap\overline P$ is a subbundle in $T_M^\C$;

ii) $P+\overline P$ is involutive.\\
Then $P$ is a strong polarization of $\omega$.
\end{propn}
\begin{proof}
The integrability of $P$ follows from the Frobenius-Nirenberg theorem.
The exactness of the Dolbeault complex for $P$ is
a theorem due to Rawnsley (\cite{Ra}, Thm. 2).
\end{proof} 

Note that in the analytic case, i.e. when $M$ is an analytic manifold and $\omega_0$ 
is an analytic form, a Lagrangean involutive analytic subbundle $P$ is 
strong. This follows from the Frobenius theorem.

\subsection{Local properties of strong polarizations}

From now on we will denote by $P$ either a polarization
of $\omega_0$ or a deformed polarization of $\omega$,
and $\OO=\OO_P$ is the sheaf of functions constant along $P$.
By $\{\cdot,\cdot\}$ we will denote the Poisson bracket
corresponding to $\omega_0$ or $\omega$.  

If $P$ is integrable, there exist, locally,
functions $a_i$, $i=1,...n$, such that $\{a_i,a_j\}=0$
and $X_{a_i}=\{a_i,\cdot\}$ form a local basis in $P$.
Denote by $\Bar d^i$ the dual basis in $P^*$.
It is clear that $\Bar d^i$ is a closed form in the
Dolbeault complex and 
\be{}\label{expd}
\Bar d=\sum_i\{a_i,\cdot\}\Bar d^i.
\ee{}

\begin{lemma}\label{lem3.1}
Suppose $P$ is a strong polarization.
Then, locally, there exist functions $f_i$ such that $\{f_i,a_j\}=\delta_{ij}$.
\end{lemma}

\begin{proof}
For any $i$, consider the closed form $\Bar d^i$. 
By Definition \ref{defn3.1}, there is a function
$f_i$ such that $\Bar d(f_i)=-\Bar d^i$. 
By (\ref{expd}), $f_i$ satisfies the lemma.
\end{proof}

Let $F_\OO=\{b\in\Cf; \{b,\OO\}\subset\OO\}$.
It follows from the Jacobi identity that
$\{F_\OO,F_\OO\}\subset F_\OO$. Hence,
$F_\OO$ is a Lie algebra acting on $\OO$ by derivations.
The kernel of this action is equal to $\OO$ itself.

Note that $F_\OO$ is a locally free $\OO$ module.
Locally, it is freely generated as an $\OO$ module 
by functions $1$, $f_1,...,f_n$, where $f_i$ are as in
Lemma \ref{lem3.1}.

Let us denote by $\Omega^1=\Omega^1_\OO$ the $\OO$-submodule of
$(T^\C)^*$ generated by the sheaf $d\OO$, i.e. $\Omega^1=\OO d\OO$. 
Let $T_\OO=F_\OO/{\OO}$. It is a locally free $\OO$ module
and is a sheaf of Lie algebras acting on $\OO$ with the local
basis $\{f_i,\cdot\}$. It is easy to see that $\Omega^1$
can be naturally identified
with $Hom_{\OO}(T_\OO,\OO)$ and $d^i=da_i$ form the dual
basis to $\{f_i,\cdot\}$. 

We form in the obvious way an analog of de Rham complex:
$(\bigwedge_{\OO_P}^\bullet\Omega^1, d)$, which we call $\OO_P$-de Rham complex.

\begin{lemma}\label{lemm3.2}
Let $P$ be a strong polarization. Then for any point of $M$ there exists a neighborhood
$U$ such that the complex $(\Gamma_U(\bigwedge_{\OO_P}^\bullet\Omega^1), d)$
is exact and gives a resolution for $\C$.
\end{lemma}

\begin{proof}
Let us choose $U$ such that the Dolbeault complex
and the usual $C^\infty$-de Rham complex $(\bigwedge^\bullet (T^\C)^*, d_{dR})$
are exact.
On the de Rham complex $(\bigwedge^\bullet (T^\C)^*, d_{dR})$
consider the Hodge filtration generated by $P^\perp\subset (T^\C)^*$.
Considering the ``stupid'' filtration on $(\bigwedge_{\OO_P}^\bullet\Omega^1, d)$
we have that the natural inclusion 
$(\bigwedge_{\OO_P}^\bullet\Omega^1, d)\to (\bigwedge^\bullet (T^\C)^*, d_{dR})$
is a filtered quasiisomorphism, i.e. induces a
quasiisomorphism on the associated graded complexes. This follows from
the fact that, since $P$ is strong,  the Dolbeault complex is exact. Now our lemma
follows from the exactness of $(\bigwedge^\bullet (T^\C)^*, d_{dR})$.
\end{proof}
 
In the following, saying that the Dolbeault or $\OO$-de Rham complex of sheaves
is exact we mean that for any point of $M$ there exists a neighborhood $U$
such that the corresponding complex of sections over $U$ is exact.

\begin{lemma}\label{lem3.2}
Let $P$ be a strong polarization.
Then $f_i$ in Lemma \ref{lem3.1} can be
chosen in such a way that $\{f_i,f_j\}=0$.
\end{lemma}

\begin{proof}
By Lemma \ref{lem3.1} one can choose $g_i$ such
that $\{g_i,a_j\}=\delta_{ij}$. For any functions
$b_1,...,b_n$ from $\OO$ the functions $f_i=g_i+b_i$
also satisfy Lemma \ref{lem3.1}.
Let us find $b_i$ such that $\{f_i,f_j\}=0$.
We find $b_i$ from the equation
$$\{g_i+b_i,g_j+b_j\}=\{g_i,g_j\}+\{g_i,b_j\}-\{g_j,b_i\}=0.$$
Note that $g=\sum_{i,j}\{g_i,g_j\}d^i\wedge d^j$ is a closed 2-form of
the $\OO_P$-de Rham complex and the equation can
be rewritten as
$$db=-g,$$
where $b=\sum_ib_id^i$.
By Lemma \ref{lemm3.2}, such $b$ exists.
\end{proof} 

Lemma \ref{lem3.2}, in particular, means that locally
there exists a polarization, $P'$, of $\omega$ complement to $P$.
This polarization is defined as annihilator of $df_1,...,df_n$.
It is obvious that locally the form $\omega$
can be written as $\omega=\sum_idf_i\wedge da_i$.
One can prove that $P'$ is also a strong polarization.

\begin{propn}\label{prop3.3}
Let $P$ be a strong polarization.
Then 

a) $P$ is a good polarization and 

b) the sequence
\be{}\label{sLa}
0\longrightarrow \OO_P\longrightarrow F_\OO\stackrel{\pi}{\longrightarrow} 
T_\OO\longrightarrow 0
\ee{}
locally splits as an exact sequence of $\OO$ modules and Lie algebras.
\end{propn}

\begin{proof}
a) follows from the existence of functions $f_i$ satisfying
Lemma \ref{lem3.2}.

b) One needs to construct, locally, a Lie algebra morphism
$s:T_\OO\to F_\OO$ such that $\pi s=id$.
Let us choose functions $a_i$ and $f_i$, $i=1,...,n$, as in Lemma \ref{lem3.2}
and assign $\{f_i,\cdot\}\mapsto  f_i$.
Since the elements $X_{f_i}=\{f_i,\cdot\}$ form a local basis in $T_\OO$,
this defines a map of $\OO_P$ modules, $s:T_\OO\to F_\OO$.
One has $[X_{f_i},X_{f_j}]=0$ and also
$\{f_i,f_j\}=0$, which proves that $s$ is a Lie algebra morphism.
\end{proof}

\begin{propn}\label{prop3.4}
Let $(M, \omega,\PP)$ and $(M,\omega',\PP')$ be two deformations of
$(M,\omega_0,P)$. Then any point of $M$ has an open neighborhood, $U$,
such that there exists a formal automorphism, $B$, of $U$ identical modulo $t$
that transforms $(U, \omega,\PP)$ to $(U, \omega',\PP')$.
\end{propn}

\begin{proof}
For any point of $m\in M$ there is an open neighborhood, $V$, 
such that on it $\omega'=\omega+td\gamma$, where
$\gamma$ is an 1-form. Let $X_\gamma$ be 
the vector field on $V$ corresponding to $\gamma$ by the isomorphism
$T_M^*\to T_M$ determined by $\omega$. Then, $\exp(tX_\gamma)$ is
a formal automorphism of $V$ which transforms $\omega$ to $\omega'$.
Hence, one may suppose in the proposition that $\omega=\omega'$ and
$\PP$ and $\PP'$ are two polarizations of the same form $\omega$.      

Let $U$ be a neighborhood of $m$ such that over $U$
the Dolbeault complexes corresponding to both $\PP$ and $\PP'$ 
are exact.
Let functions $a_i$, $a'_i$ on $U$ be such that
$da_i$ and $da'_i$ form local basises in $\PP^\perp$ and $\PP^{'\perp}$,
respectively. Suppose $\PP=\PP'$ mod $t^k$. 
Then, we may assume that $a_i=a'_i$ mod $t^k$, i.e. 
$a'_i=a_i+t^kb_i$ for some functions $b_i$ on $U$.
Let $\{\cdot,\cdot\}$ denote the Poisson bracket inverse to $\omega$.
One has $\{a_i,a_j\}=\{a'_i,a'_j\}=0$, hence $\{a_i,b_j\}+\{b_i,a_j\}=0$.
It follows from the exactness of Dolbeault complex that there exists a function $g$
on $U$ such that $\{a_i,g\}=b_i$ for all $i$.
Let $X_g$ be the Hamiltonian vector field on $U$ corresponding to $g$.
Then, the automorphism $\exp(tX_g)$ of $U$ leaves $\omega$ on the place
and transforms $\PP$ to a polarization, $\PP''$, such that 
$\PP''=\PP'$ mod $t^{k+1}$. Proceeding by iteration proves the proposition.  
\end{proof}

This proposition shows that any deformation
$(M, \omega,\PP)$ of $(M, \omega_0,P)$ is, locally,
isomorphic to the trivial deformation $(M, \omega_0,P[[t]])$.

\subsection{Local properties of polarized star-products}

Let $(M,\omega_0)$ be a symplectic manifold with a strong 
polarization $P$. Let $\omega$ be a symplectic deformation of $\omega_0$
and $\PP$ a polarization of $\omega$ which is a deformation of $P$.
As is proven in the previous section, $\PP$ is strong and thus a good
polarization. Let $\OO$ be the sheaf of functions constant along $\PP$  
and $(\AA,\OO)$ a polarized star-product constructed in
Section 2.

Our immediate goal is to prove analogs of Lemmas \ref{lem3.1} and \ref{lem3.2}
for the commutator $\ft[\cdot,\cdot]$ in $\AA$.

\begin{lemma}\label{lem4.1}
Let $\PP$ be a strong polarization of $\omega$.
Let $a_1,...,a_n\in \OO$ be functions on an open set $U\subset M$ such that 
$da_1,...,da_n$ form a local basis in $\PP^\perp$.
Then, locally, there exist functions $\cf_1,...,\cf_n$ such that
$\ft[\cf_i,a_j]=\delta_{ij}$.
\end{lemma}  

\begin{proof}
The operators $\ft[a_i,\cdot]$ are pairwise commuting derivations
of $\AA$ restricted to an open set $U$. 
Let $L$ be the free $\CM$ module over $U$ spanned on $[a_i,\cdot]$.
Denote by $\hat d^i\in L^*$ the basis dual to $\ft[a_i,\cdot]$.
Let us define in the obvious way the complex
$(\bigwedge^\bullet_{\CM} L^*, \hat d)$, where
\be{}\label{expd1}
\hat d=\sum_i[a_i,\cdot]\hat d^i.
\ee{}
Since $\ft[a_i,\cdot]=\{a_i,\cdot\}+o(t)$ and $\rank L=\rank P$, 
this complex is a deformation of the Dolbeault complex 
$(\bigwedge^{\bullet}_{\CM} P^*, \Bar d)$ over $U$ with $\Bar d$ defined
by (\ref{expd}).
Since $\PP$ is strong, the Dolbeault complex is exact, therefore
the complex $(\bigwedge^\bullet_{\CM} L^*, \hat d)$ is exact, too.

Now we proceed as in the proof of Lemma \ref{lem3.1}.
Namely, for any $i$, consider the closed form $\hat d^i$. 
Since the complex $(\bigwedge^\bullet_{\CM} L^*, \hat d)$
is exact,
there is a function
$\cf_i$ such that $\hat d(\cf_i)=-\hat d^i$. 
By (\ref{expd1}), $f_i$ satisfies the lemma.
\end{proof}

Let
\be{}
 F(\AA)=\{b\in\CM; [b,\OO]\subset\OO\}.
\ee{}
It is clear that $F(\AA)$ is a sheaf of Lie algebras on $M$ with
the bracket $\ft[\cdot,\cdot]$.
It is a locally free $\OO$ module with
the local basis consisting of functions $1$, $\cf_j$.
$\OO$ sits in $F(\AA)$ as a commutative Lie subalgebra
The sheaf of Lie algebras $F(\AA)/\OO$ is also a locally free $\OO$-module
isomorphic to the sheaf
$T_\OO=F_\OO/\OO$ from the previous section:
locally, the isomorphism is given by $\ft[\cf,\cdot]\mapsto \{f_i,\cdot\}$.

\begin{lemma}\label{lem4.2}
By hypothesis of Lemma \ref{lem4.1},
$\cf_i$ can be chosen in such a way that
$\ft[\cf_i,\cf_j]=0$.
\end{lemma}
\begin{proof} The same as of Lemma \ref{lem3.2}.
\end{proof} 

Lemma \ref{lem4.2}, in particular, means that locally
there exists a ``complement'' commutative subalgebra in $\AA$.
This subalgebra is generated by functions $\cf_i$.
But that ``complement'' subalgebra may {\em not} be an algebra
of functions constant on a polarization.

\begin{propn}\label{prop4.3}
Let $\PP$ be a strong polarization of $\omega$ and $(\AA,\OO)$ 
the corresponding polarized quantization.
Then the sequence
\be{}\label{sLa1}
0\longrightarrow \OO\longrightarrow F(\AA)\stackrel{\pi}{\longrightarrow} 
T_\OO\longrightarrow 0
\ee{}
locally splits as an exact sequence of $\OO$ modules and Lie algebras.
\end{propn}

\begin{proof}
The same as of Proposition \ref{prop3.3}.
\end{proof}

\begin{propn}\label{prop4.4}
Let $(M, \omega,\PP)$, $(M,\omega',\PP')$ be two deformations of
$(M,\omega_0,P)$ and $(\AA,\OO)$, $(\AA',\OO')$
corresponding polarized star-products.
Then any point of $M$ has an open neighborhood, $U$,
such that over $U$ these star-products are equivalent.
Moreover, if $\OO=\OO'$ then
the equivalence can be chosen to be identity on $\OO$.
\end{propn}

\begin{proof}
By Proposition \ref{prop3.4} one can suppose that
$\omega=\omega'=\omega_0$ and $\PP=\PP'=P[[t]]$. 
So, $\OO=\OO'=\OO_0[[t]]$, where $\OO_0$ is the algebra of 
functions constant along $P$.

Let $m\in M$ and $U$ be a sufficiently small neighborhood of $m$.
Let $a_i$ be functions on $U$ such that $da_i$ form a local basis in $P^\perp$
and $f_i$ be functions such that $\{f_i,a_j\}=\delta_{ij}$, where
$\{\cdot,\cdot\}$ is the Poisson bracket inverse to $\omega_0$.
Let $\mu=\sum_{i\geq 0}t^i\mu_i$ and $\mu'=\sum_{i\geq 0}t^i\mu'_i$
be multiplications in $\CM$ over $U$ corresponding to $\AA$ and $\AA'$ 
and consisting of bidifferential operators. So, $\mu_0(a,b)=\mu'_0(a,b)=ab$ and 
one may suppose that $\mu_1(a,b)=\mu'_1(a,b)=\frac{1}{2}\{a,b\}$.

Suppose that there exists a differential operator $B^{(n)}=1+\cdots +t^n B_n$
which transforms $\mu$ to $\mu'$ modulo $t^{n+1}$ and such
that $B_k$, $k=1,...,n$, take $\OO$ to zero.
This means, in particular, that we may suppose that $\mu=\mu'$ modulo $t^{n+1}$, so that
$$\mu=\mu_0+\cdots t^n\mu_n+\mu_{n+1}+\cdots,$$
$$\mu'=\mu_0+\cdots t^n\mu_n+\mu'_{n+1}+\cdots.$$
We are going to prove that there exists a differential operator of
the form $B=1+t^nX+t^{n+1}B_{n+1}$, where $X$ is a vector field on $U$, 
which transforms $\mu'$ to $\mu$ modulo $t^{n+2}$  
and such that $X(\OO)=B_{n+1}(\OO)=0$.

It is easy to check that $c=\mu'_{n+1}-\mu_{n+1}$ is a Hochschild cocycle
in the algebra $C^\infty_U$.
Hence, $\nu(a,b)=c(a,b)-c(b,a)$ is a bivector field
and there is the representation: $\mu'_{n+1}=\mu_{n+1}+\nu+\delta_H(B^\prime_{n+1})$,
where $B_{n+1}$ is a differential operator and $\delta_H$ the Hochschild
differential.

It is also easy to check
that $\[\mu_1,\nu\]=0$, where $\[\cdot,\cdot\]$ denotes the Schouten
bracket of polyvector fields. So, there exists a vector field $Y$ such that
$\[\mu_1,Y\]=\nu$. 
This means that for any $f,g\in C^\infty_U$
$$Y(\mu_1(f,g))-\mu_1(Yf,g)-\mu_1(f,Yg)=\nu(f,g).$$
Recall now that $\mu_1(a,b)=\nu(a,b)=0$ for
$a,b\in \OO$. It follows that for any $a,b\in\OO$
$$\mu_1(Ya,b)+\mu_1(a,Yb)=0.$$
In particular, we have $\mu_1(Ya_i,a_j)+\mu_1(a_i,Ya_j)=0$
for any $a_i,a_j$ of our basis.
It follows from the exactness of Dolbeault complex
that there exists a function $b\in C^\infty_U$ such that
$Ya_i=\mu_1(b,a_i)$ for all $i$. Denote $X_b=\{b,\cdot\}$, the Hamiltonian
vector field corresponding to $b$. Put $X=Y-X_b$. Then,  
since $\[\mu_1,X_b\]=0$, we get $\[\mu_1,X\]=\nu$.
Since $Ya_i=X_ba_i$ for all $i$, we get $X(\OO)=0$. 

If we transform the multiplication $\mu'$ by the operator $1+t^nX$,
we obtain the multiplication $\mu''$
of the form $\mu''=\mu_0+\cdots t^n\mu_n+\mu''_{n+1}+\cdots$,
where $\mu''_{n+1}=\mu_{n+1}+\delta_HB^\prime_{n+1}$.

Since $\mu''_{n+1}(a,b)=\mu_{n+1}(a,b)=0$ for all $a,b\in\OO$,
we have that $B^\prime_{n+1}$ is a derivation from $\OO_U$ to $C^\infty_U$.
This derivation can be extended to a derivation of $C^\infty_U$.
Indeed, put $B^\prime_{n+1}(a_i)=g_i$. Then the operator
$D=g_i\{f_i,\cdot\}$ is such an extension. Now put 
$B_{n+1}=B^\prime_{n+1}-D$. One has $\delta_HB_{n+1}=\delta_HB^\prime_{n+1}$.
It is easy to see that the operator $1+t^{n+1}B_{n+1}$
transforms $\mu''$ to $\mu$ modulo $t^{n+2}$.
It follows that the operator
$(1+t^{n+1}B_{n+1})(1+t^nX)$ transforms $\mu'$ to $\mu$ modulo $t^{n+2}$.
\end{proof}

Let $(M,\omega_0,P)$ be a polarized symplectic manifold with
$P$ a strong polarization. Let $U$ be a sufficiently small neighborhood of a point of $M$.
Let $a_i$ be functions on $U$ such that $da_i$ form a local basis in $P^\perp$
and $f_i$ be functions such that $\{f_i,a_j\}=\delta_{ij}$, where
$\{\cdot,\cdot\}$ is the Poisson bracket inverse to $\omega_0$. 
Then the bivector field 
$$\pi=\sum_i X_{a_i}\wedge X_{f_i}$$
represent this Poisson bracket.
Since all $X_{a_i}$ and $Y_{f_i}$ pairwise commute,
one can construct the Moyal star-product of the form
\be{}\label{Moyal}
f\star g=\mu_0\exp(\frac{t}{2}X_{a_i}\wedge X_{f_i})(f\ot g),
\ee{}
where $\mu_0$ is the usual multiplication.
This star-product is obviously polarized, and we call it
the {\em polarized Moyal star-product}.

\begin{corollary}\label{cor3.1}
Any polarized star-product is
locally equivalent to a polarized Moyal star-product.
\end{corollary}

\begin{proof} Follows from the previous proposition.
\end{proof} 

\section{A formula for polarized quantization}

\subsection{$\OO$-extension of $T_\OO$}
The sequences (\ref{sLa}) and (\ref{sLa1}) are examples of an $\OO$-extension of $T_\OO$.
In the following we suppose that $\OO$ is the sheaf of functions constant along
a strong polarization and $T_\OO=Hom_\OO(\Omega^1_\OO,\OO)$.
As we have seen, $T_\OO$ is a locally free $\OO$ module of derivations of $\OO$. 

Following \cite{BK}, \cite{BB}, 
we give the following

\begin{definition}\label{defLA}
An {\em $\OO$-extension of $T_\OO$} is a locally free
sheaf $\widetilde T$ of $\OO$-modules equipped with a structure of a
sheaf of Lie algebras over $\Ct$ (with the bilinear operation denoted,
as usual, by $[\cdot,\cdot]$), a section $\cc$ of the center of $\widetilde T$,
and a surjective $\OO$-linear map
$\sigma : \widetilde T \to T_\OO$ (the anchor map) which is a Lie algebra
homomorphism, whose kernel (which is a sheaf of $\OO$-Lie algebras) is
isomorphic to $\OO$ (the latter equipped with the trivial Lie bracket).
In addition the Leibnitz rule holds: for $f\in\OO$, $\tau_1,\tau_2\in\widetilde T$,
$[\tau_1,f\tau_2]=f[\tau_1,\tau_2]+\sigma(\tau_1)(f)\tau_2$.
Thus, there is an exact sequence (of $\OO$-modules and Lie algebras)
\be{}\label{PLA}
0\longrightarrow\OO\stackrel{i}{\longrightarrow}\widetilde T\stackrel{\sigma}{\longrightarrow} 
T_\OO\longrightarrow 0 \ ,
\ee{}
where $i(f)=f\cdot\cc$, $f\in\OO$.
\end{definition}

Just as in Proposition \ref{prop3.3}, one shows that the extension \eqref{PLA}
admits local splittings which are $\OO$-linear Lie algebra homomorphisms.
Such a splitting is called a {\em flat connection on $\widetilde T$}.

\subsection{Characteristic class of an $\OO$-extension of $T_\OO$}
To each $\widetilde T$ one associates the cohomology class
$c(\widetilde T)\in H^1(M;\Omega_\OO^{1,cl})= H^2(M,\Ct)$ as follows.

Let $\{U_\alpha\}$ be a sufficiently fine open covering of $M$, so that
there are flat connections $s_\alpha$ on $\widetilde T\vert_{U_\alpha}$.
If $U_{\alpha\beta}\stackrel{def}{=}U_\alpha\cap U_\beta\neq\emptyset$,
the difference $s_\alpha - s_\beta$ gives rise to the section $c_{\alpha,\beta}$
(defined over $U_{\alpha\beta}$) of $\Omega_\OO^{1,cl}= \shHom_\OO(T_\OO,\OO)$
by the formula $i(c_{\alpha,\beta}(\xi))= s_\alpha(\xi)-s_\beta(\xi)$.
The collection $\{c_{\alpha,\beta}\}$ is a degree one \^ Cech cochain with
coefficients in the sheaf $\Omega_\OO^{1,cl}$. It is, in fact, a cocycle,
whose cohomology class, denoted $c(\widetilde T)$ is independent of the choice
of local flat connections.

The above construction recovers the Chern class of a complex line bundle.
Namely, suppose that $L$ is a locally free sheaf of $\OO$-modules of rank
one. Let $\widetilde T_{\cal L}$ denote the sheaf of differential operators on $L$
of order (at most) one. The sheaf $\widetilde T_{\cal L}$ is equipped with the
{\em left} $\OO$-module structure, the Lie bracket given by the commutator,
the central section -- the identity operator (so that the map $i$ is simply the
inclusion of operators of order zero). The principal symbol map serves as the anchor
map $\sigma$. These data  exhibit $\widetilde T_{\cal L}$ as an $\OO$-extension of $T_\OO$.

On the other hand, the sheaf ${\cal L}\otimes_\OO\Cf$ is the sheaf of $C^\infty$
sections of a complex line bundle $L$. It is easy to see that the class
$c(\widetilde T_{\cal L})$ coincides with the the first Chern class of $L$.

For example, since $\det\PP^\perp = \Omega^n_\OO\otimes_\OO\Cf$, we find that
$c_1(\PP^\perp)=c(\widetilde T_{\Omega^n_\OO})$.

For a PDQ (polarized deformation quantization) $(\AA,\OO)$ we set $\cl(\AA,\OO) = c(F(\AA))$.

\newcommand{\op}{\dagger}
\subsection{The opposite PDQ}
For a ring $R$ we denote by $R^{op}$ the opposite ring. That is, $R^{op}$ is a ring,
there is a bijection $(\bullet)^\op :R\rightarrow R^{op}$ such that $1^\op = 1$ and
$r_1^\op r_2^\op = (r_2r_1)^\op$. We will usually identify $R$ and $R^{op}$ using the
bijection.

If $\AA$ is a deformation quantization of $(M,\omega_0)$, then, clearly, $\AA^{op}$
is a deformation quantization of $(M,-\omega_0)$. If $(\AA,\OO)$ is a PDQ, then
$(\AA^{op},\OO)$ is a PDQ (since $(\bullet)^\op$ restricts to a ring isomorphism
on $\OO$). It is clear that $F(\AA)$ and $F(\AA^{op})$ coincide as subsheaves of
$\AA$. However, the respective $\OO$-module structures, symbol maps, and Lie brackets
are different. In fact, one has:
\begin{itemize}
\item $f\tau^\op = (f\tau -t\sigma(\tau)(f))^\op = (f\tau)^\op -t\sigma(\tau)(f)$
\item $\sigma(\tau^\op)=-\sigma(\tau)$
\item $[\tau_1^\op,\tau_2^\op]=-[\tau_1,\tau_2]^\op$.
\end{itemize}

\begin{propn}\label{pr4.1}
Suppose that $(\AA,\OO)$ is a PDQ. Then,
\be{*}
\cl(\AA^{op},\OO) = c_1(\PP^\perp) - \cl(\AA,\OO)
\ee{*}
\end{propn}
\begin{proof}
Choose an open covering $\{U_\alpha\}$ of $M$ such that, for each open set
$U_\alpha$ there exist flat connections $s_\alpha$ on $F(\AA)\vert_{U_\alpha}$ and
$s_\alpha^\Omega$ on $\widetilde T_{\Omega^n_\OO}\vert_{U_\alpha}$.

Let $\ell$ denote the splitting of the extension
\be{*}
0\longrightarrow\OO\longrightarrow\widetilde T_{\Omega^n_\OO}\longrightarrow 
T_\OO\stackrel{\sigma}{\longrightarrow} 0 \ 
\ee{*}
given by the {\em negative} of the Lie derivative. This splitting is a Lie algebra
homomorphism, but is {\em not} $\OO$-linear.
Instead, for $f\in\OO$ and $\xi\in T$, $\ell(f\xi)= f\ell(\xi)+\xi(f)$.
We denote the restriction of this splitting to $U_\alpha$ by the same letter.

Since both $\ell$ and $s_\alpha^\Omega$ are splittings of the same sequence,
their difference $\psi_\alpha = s_\alpha^\Omega - \ell$ is a map
$\psi_\alpha : T \to \OO$ which is {\em not} $\OO$-linear, but satisfies
$\psi_\alpha(f\xi)= f\psi_\alpha(\xi)+\xi(f)$. Note that $\psi_\alpha -\psi_\beta
= s_\alpha^\Omega - s_\beta^\Omega$.

Let $s_\alpha^{op}(\xi)=t\psi_\alpha(\xi) - s_\alpha(\xi)^\op$. We claim that
$s_\alpha^{op}$ is a (locally defined) flat connection on $F(\AA^{op})$. To this
end we check various properties.

The calculation
\be{*}
\sigma(s_\alpha^{op}(\xi))=\sigma(t\psi_\alpha(\xi) - s_\alpha(\xi)^\op)=
\sigma(- s_\alpha(\xi)^\op)=\sigma(s_\alpha(\xi))=\xi
\ee{*}
shows that $s_\alpha^{op}$ is a (local) splitting. The calculation
\begin{multline*}
s_\alpha^{op}(f\xi)=t\psi_\alpha(f\xi) - s_\alpha(f\xi)^\op=
\psi_\alpha(f\xi) - (fs_\alpha(\xi))^\op=
tf\psi_\alpha(\xi)+t\xi(f)-(fs_\alpha(\xi)^\op+t\xi(f))= \\
tf\psi_\alpha(\xi)-fs_\alpha(\xi)^\op=fs_\alpha^{op}(\xi)
\end{multline*}
shows that this splitting is $\OO$-linear.  

 One verifies other properties
similarly.

Hence one can use the locally defined flat connections $\{s_\alpha^{op}\}$
to calculate the characteristic class of $F(\AA^{op})$. The desired formula
follows from
\be{*}
s_\alpha^{op} - s_\beta^{op} = (\psi_\alpha(\xi) - s_\alpha(\xi)) -
(\psi_\beta(\xi) - s_\beta(\xi) = (s_\alpha^\Omega - s_\beta^\Omega) -
(s_\alpha - s_\beta)
\ee{*}
\end{proof}

\subsection{$\cl(\AA,\OO)$ and other classes}

Let $\AA_1$ and $\AA_2$ be two deformation quantizations on $(M,\omega_0)$.
Deligne, \cite{De},  related to them cohomology classes
\be{*}
[\AA_2:\AA_1]\in H^2(M,\Ct) \qquad \mbox{\rm and}\\
\obs(\AA_1)\in \ft[\omega_0]+ H^2(M,\Ct). 
\ee{*}
He proved that
\be{}\label{Del}
[\AA_2:\AA_1]=\theta(\AA_2)-\theta(\AA_1) \qquad\mbox{\rm and}\\
\obs(\AA_1)=-t\frac{d}{dt}\theta(\AA_1). \label{Del1}
\ee{}

\begin{lemma}\label{lem5.1}
Let  $(\AA_1,\OO_1)$ and  $(\AA_2,\OO_2)$ be
two polarized quantizations of $(M,\omega_0,P)$.
Then $[\AA_2:\AA_1]=\ft\cl(\AA_2,\OO_2)-\ft\cl(\AA_1,\OO_1)$.
\end{lemma}

\begin{proof}
By Proposition \ref{prop4.4}, 
there is an open covering  $\{U_\alpha\}$ of $M$ such that
there exist isomorphisms of polarized quantizations
\be{*}
B_\alpha:(\AA_1,\OO_1)\mid_{U_\alpha}\to (\AA_2,\OO_2)\mid_{U_\alpha}.
\ee{*}
These isomorphisms induce isomorphisms
$F(\AA_1)\to F(\AA_2)$ over each $U_\alpha$.
We may suppose that
both $F(\AA_1)$ and $F(\AA_2)$  admit
on each $U_\alpha$ flat connections
$s_{i\alpha}:T_{\OO_i}\to F(\AA_i)$, $i=1,2$,
such that $B_\alpha(s_{1\alpha})=s_{2\alpha}$. 
Since $F(\AA_1)=F(\AA_2)$ mod $t$,
we may choose the connections in such a way that
$s_{1\alpha}=s_{2\alpha}$ mod $t$.

Let $f_{i\alpha\beta}\in\OO_i$ be functions on $U_\alpha\cap U_\beta$ such that
$df_{i\alpha\beta}=s_{i\beta}-s_{i\alpha}$.
Since $s_{1\alpha}=s_{2\alpha}$ mod $t$,
we may choose $f_{i\alpha\beta}$ such that 
$f_{1\alpha\beta}=f_{2\alpha\beta}$ mod $t$.
On $U_\alpha\cap U_\beta$ one has
$B^{-1}_\alpha B_\beta=\exp(\ft ad(f_{2\alpha\beta}-f_{1\alpha\beta}))$.
Since $f_{2\alpha\beta}-f_{1\alpha\beta}\in t\OO_i$,
$\exp(\ft ad(f_{2\alpha\beta}-f_{1\alpha\beta}))$ is an
automorphism of $\AA_1$ on $U_\alpha\cap U_\beta$.
According to Deligne's definition,
$[\AA_2:\AA_1]$ is equal to $\ft\check\partial (f_{2\alpha\beta}-f_{1\alpha\beta})$.
\end{proof}

\begin{lemma}\label{lem5.2}
Let $(\AA,\OO)$ be a polarized quantization of $(M,\omega_0,P)$
corresponding to a deformation $(M,\omega,\PP)$.
Then, $\obs(\AA)=-t\frac{d}{dt}(\ft\cl(\AA,\OO))$.
\end{lemma}

\begin{proof}
First, we assume that $(\AA,\OO)$ corresponds to the trivial deformation
$(M,\omega_0,P[[t]])$.

By Proposition \ref{prop4.4}, 
there is an open covering  $\{U_\alpha\}$ of $M$ such that
$(\AA,\OO)\mid_{U_\alpha}$ are isomorphic to Moyal polarized quantizations.
Let $D_\alpha$ be $\C$ linear derivations of $\AA\mid_{U_\alpha}$
such that for $a\in\OO$ one has $D_\alpha a=t\dt a$.
Such $D_\alpha$ exist because they obviously exist for the Moyal star-products.
Then, $D_\alpha F(\AA)\subset F(\AA)$.

Let the covering $\{U_\alpha\}$ be fine enough and 
on each $U_\alpha$ there exists a flat connection
$s_{\alpha}:T_{\OO}\to F(\AA)$. 
Then, on $U_\alpha\cap U_\beta$ one has
$s_\beta-s_\alpha=df_{\alpha\beta}$ for $f_{\alpha\beta}\in\OO$. 
Let $a_{\alpha 1},...,a_{\alpha n}$ be functions not depending on $t$
such that 
$da_{\alpha 1},...,da_{\alpha n}$ 
form a local basis in $P^\perp$,
and $\partial_{\alpha 1},...,\partial_{\alpha n}$
be the dual basis in $T_\OO$. Put $f_{\alpha i}=s_{\alpha i}(\partial_{\alpha i})$.

One has $[f_{\alpha i},a_{\alpha j}]=t\delta_{ij}$. Hence,
since $a_{\alpha i}$ does not depend of $t$,
$[D_\alpha f_{\alpha i},a_{\alpha j}]=t\delta_{ij}$.
This means that 
$D_\alpha f_{\alpha i}=f_{\alpha i} +b_{\alpha i}$ for some
$b_{\alpha i}\in\OO$.
It follows from $[f_{\alpha i},f_{\alpha j}]=0$ that
$[f_{\alpha i},b_{\alpha j}]+[b_{\alpha i},f_{\alpha j}]=0$.
It follows from the exactness of $\OO$-de Rham complex
that there exist $b_\alpha\in\OO$ such that
$b_{\alpha i}=\ft[b_\alpha,f_{\alpha i}]$.
So, $D_\alpha f_{\alpha i}=f_{\alpha i}+\ft[b_\alpha,f_{\alpha i}]$.
Taking $D_\alpha-\ft\ad(b_\alpha)$ instead $D_\alpha$
we obtain that the new $D_\alpha$ satisfy
$D_\alpha f_{\alpha i}=f_{\alpha i}$.
This equations may be rewritten in the form
\be{}\label{goodr}
D_\alpha s_{\alpha }=s_{\alpha }.
\ee{}
Further, one has
\be{}\label{vr}
s_\alpha=s_\beta-[\ft f_{\alpha\beta},s_\beta].
\ee{}
Using (\ref{goodr}) and (\ref{vr}), we have
$$(D_\beta-D_\alpha)s_\alpha=
D_\beta(s_\beta-[\ft f_{\alpha\beta},s_\beta])-D_\alpha s_\alpha=
-[D_\beta(\ft f_{\alpha\beta}),s_\beta]=[-t\dt(\ft f_{\alpha\beta}),s_\beta].$$
Since derivations of $\AA$ are completely defined by their
values on $F(\AA)$, it follows that 
$D_\beta-D_\alpha=\ad(-t\dt(\ft f_{\alpha\beta})).$
According to Deligne's definition, this means that
$(-t\dt(\ft f_{\alpha\beta}))$ represents $\obs(\AA)$.

Now, let $(M,\omega,\PP)$ be arbitrary.
Denote by $\cl$ the class of a quantization corresponding to $(M,\omega,\PP)$
and by $\cl_0$ the class of a quantization corresponding to 
the trivial deformation $(M,\omega_0,P[[t]])$. The same meaning has
the notation $\theta$, $\theta_0$ and $\obs$, $\obs_0$.

By Lemma \ref{lem5.1} and (\ref{Del}) we have
$\ft\cl_0-\ft\cl=\theta_0-\theta$, hence
$$\theta_0-\ft\cl_0=\theta-\ft\cl.$$
Applying $-t\dt$ to the left hand side of this equation, we obtain zero
(by (\ref{Del1}) and what we have just proved).
Hence, applying $t\dt$ to the right hand side and using (\ref{Del1}), we have 
$$\obs=-t\dt\theta=-t\dt\cl,$$
which proves the proposition.
\end{proof}

It follows from Lemma \ref{lem5.1} that the difference
$\theta(\AA)-\ft\cl(\AA,\OO)$ does not depend on the
quantization $(\AA,\OO)$ of $(M,\omega_0,P)$,
while Lemma \ref{lem5.2} shows that this difference
does not depend on $t$, i.e. 
\be{}\label{prefor}
\theta(\AA)=\ft\cl(\AA,\OO)-c,
\ee{}
where $c\in\C$.

\subsection{Relation between $\theta(\AA)$ and $\cl(\AA,\OO)$ }
In this Section we calculate the constant $c$ in (\ref{prefor}).

Let $\AA$ be a deformation quantization of $(M,\omega_0)$.
Let $\AA^{op}$ denote the opposite algebra
and $\AA^\sigma$ denote the algebra obtained by
the automorphism of $\Ct$ taking $t\mapsto -t$. 
It is clear that both $\AA^{op}$ and $\AA^\sigma$ are
quantizations of $(M,-\omega_0)$.

\begin{lemma}\label{lem5.3}
Let $\AA$ be a deformation quantization on $(M,\omega_0)$.
Then

a)  $\theta(\AA^{op})=-\theta(\AA)$;

b) $\theta(\AA^\sigma)=\theta^\sigma(\AA)$.
\end{lemma}

\begin{proof} Follows from the Fedosov construction, see also \cite{De}.
\end{proof}

\begin{thm}\label{thm5.1}
Let $(\AA,\OO)$ be a polarized quantization of $(M,\omega_0,P)$.
Then
\be{}\label{Mf}
\theta(\AA)=\ft\cl(\AA,\OO)-\frac{1}{2}c_1(P).
\ee{}
\end{thm}

\begin{proof}
By Lemma \ref{lem5.1}  we have the equalities
$$
\theta(\AA^{op}) - \theta(\AA^\sigma) = [\AA^{op}:\AA^\sigma] =
\cl(\AA^{op},\OO) - \cl((\AA,\OO)^\sigma)\ .
$$
By Lemma \ref{lem5.3} and Proposition \ref{pr4.1} 
we have $\theta(\AA^{op}) = -\theta(\AA)$ and
$\cl(\AA^{op},\OO) = -\cl(\AA,\OO) + c_1(\Omega^1_\OO)$. 
Combining these with the preceeding
identity we obtain
$$
(\cl(\AA,\OO) - \theta(\AA)) + (\cl((\AA,\OO)^\sigma) - \theta(\AA^\sigma)) =
c_1(\Omega^1_\OO)\ .
$$
By (\ref{prefor}) both differences in the left hand side of the last identity are equal to the
differences of the respective constant terms, hence coincide (since the constant terms are not 
affected by $\sigma$). Thus the above expression says
$$
2(\cl(\AA,\OO) - \theta(\AA)) = c_1(\Omega^1_\OO)\ .
$$
The theorem follows from the fact that $c_1(\Omega^1_\OO)=c_1(P)$.
\end{proof}  

Note that the formula (\ref{Mf}) relates to the formula (\ref{Wc}), which
shows that $\cl(\AA,\OO)=\ft[\Omega_D]$, where
$\Omega_D$ is the Wick curvature appearing in the Fedosov procedure
of constructing $(\AA,\OO)$.

\bigskip

{\em e-mail:} bressler@ihes.fr

{\em e-mail:} donin@macs.biu.ac.il

\end{document}